\documentstyle[11pt]{article}
\title{\bf Invariant metrics on $G$-spaces}
\author{Bogus\l aw Hajduk, Rafa\l\ Walczak (Wroc\l aw)}
\catcode`\@=11

\font\tenmsx=msam10
\font\sevenmsx=msam7
\font\fivemsx=msam5
\font\tenmsy=msbm10
\font\sevenmsy=msbm7
\font\fivemsy=msbm5
\newfam\msxfam
\newfam\msyfam
\textfont\msxfam=\tenmsx  \scriptfont\msxfam=\sevenmsx
  \scriptscriptfont\msxfam=\fivemsx
\textfont\msyfam=\tenmsy  \scriptfont\msyfam=\sevenmsy
  \scriptscriptfont\msyfam=\fivemsy

\def\hexnumber@#1{\ifcase#1 0\or1\or2\or3\or4\or5\or6\or7\or8\or9\or
        A\or B\or C\or D\or E\or F\fi }

\font\teneuf=eufm10
\font\seveneuf=eufm7
\font\fiveeuf=eufm5
\newfam\euffam
\textfont\euffam=\teneuf
\scriptfont\euffam=\seveneuf
\scriptscriptfont\euffam=\fiveeuf
\def\frak{\ifmmode\let\next\frak@\else
 \def\next{\Err@{Use \string\frak\space only in math mode}}\fi\next}
\def\goth{\relaxnext@\ifmmode\let\next\frak@\else
 \def\next{\Err@{Use \string\goth\space only in math mode}}\fi\next}
\def\frak@#1{{\frak@@{#1}}}
\def\frak@@#1{\fam\euffam#1}

\edef\msx@{\hexnumber@\msxfam}
\edef\msy@{\hexnumber@\msyfam}

\mathchardef\boxdot="2\msx@00
\mathchardef\boxplus="2\msx@01
\mathchardef\boxtimes="2\msx@02
\mathchardef\square="0\msx@03
\mathchardef\blacksquare="0\msx@04
\mathchardef\centerdot="2\msx@05
\mathchardef\lozenge="0\msx@06
\mathchardef\blacklozenge="0\msx@07
\mathchardef\circlearrowright="3\msx@08
\mathchardef\circlearrowleft="3\msx@09
\mathchardef\rightleftharpoons="3\msx@0A
\mathchardef\leftrightharpoons="3\msx@0B
\mathchardef\boxminus="2\msx@0C
\mathchardef\Vdash="3\msx@0D
\mathchardef\Vvdash="3\msx@0E
\mathchardef\vDash="3\msx@0F
\mathchardef\twoheadrightarrow="3\msx@10
\mathchardef\twoheadleftarrow="3\msx@11
\mathchardef\leftleftarrows="3\msx@12
\mathchardef\rightrightarrows="3\msx@13
\mathchardef\upuparrows="3\msx@14
\mathchardef\downdownarrows="3\msx@15
\mathchardef\upharpoonright="3\msx@16

\mathchardef\downharpoonright="3\msx@17
\mathchardef\upharpoonleft="3\msx@18
\mathchardef\downharpoonleft="3\msx@19
\mathchardef\rightarrowtail="3\msx@1A
\mathchardef\leftarrowtail="3\msx@1B
\mathchardef\leftrightarrows="3\msx@1C
\mathchardef\rightleftarrows="3\msx@1D
\mathchardef\Lsh="3\msx@1E
\mathchardef\Rsh="3\msx@1F
\mathchardef\rightsquigarrow="3\msx@20
\mathchardef\leftrightsquigarrow="3\msx@21
\mathchardef\looparrowleft="3\msx@22
\mathchardef\looparrowright="3\msx@23
\mathchardef\circeq="3\msx@24
\mathchardef\succsim="3\msx@25
\mathchardef\gtrsim="3\msx@26
\mathchardef\gtrapprox="3\msx@27
\mathchardef\multimap="3\msx@28
\mathchardef\therefore="3\msx@29
\mathchardef\because="3\msx@2A
\mathchardef\doteqdot="3\msx@2B

\mathchardef\triangleq="3\msx@2C
\mathchardef\precsim="3\msx@2D
\mathchardef\lesssim="3\msx@2E
\mathchardef\lessapprox="3\msx@2F
\mathchardef\eqslantless="3\msx@30
\mathchardef\eqslantgtr="3\msx@31
\mathchardef\curlyeqprec="3\msx@32
\mathchardef\curlyeqsucc="3\msx@33
\mathchardef\preccurlyeq="3\msx@34
\mathchardef\leqq="3\msx@35
\mathchardef\leqslant="3\msx@36
\mathchardef\lessgtr="3\msx@37
\mathchardef\backprime="0\msx@38
\mathchardef\risingdotseq="3\msx@3A
\mathchardef\fallingdotseq="3\msx@3B
\mathchardef\succcurlyeq="3\msx@3C
\mathchardef\geqq="3\msx@
\mathchardef\geqslant="3\msx@3E
\mathchardef\gtrless="3\msx@3F
\mathchardef\sqsubset="3\msx@40
\mathchardef\sqsupset="3\msx@41
\mathchardef\vartriangleright="3\msx@42
\mathchardef\vartriangleleft="3\msx@43
\mathchardef\trianglerighteq="3\msx@44
\mathchardef\trianglelefteq="3\msx@45
\mathchardef\bigstar="0\msx@46
\mathchardef\between="3\msx@47
\mathchardef\blacktriangledown="0\msx@48
\mathchardef\blacktriangleright="3\msx@49
\mathchardef\blacktriangleleft="3\msx@4A
\mathchardef\vartriangle="0\msx@4D
\mathchardef\blacktriangle="0\msx@4E
\mathchardef\triangledown="0\msx@4F
\mathchardef\eqcirc="3\msx@50
\mathchardef\lesseqgtr="3\msx@51
\mathchardef\gtreqless="3\msx@52
\mathchardef\lesseqqgtr="3\msx@53
\mathchardef\gtreqqless="3\msx@54
\mathchardef\Rrightarrow="3\msx@56
\mathchardef\Lleftarrow="3\msx@57
\mathchardef\veebar="2\msx@59
\mathchardef\barwedge="2\msx@5A
\mathchardef\doublebarwedge="2\msx@5B
\mathchardef\angle="0\msx@5C
\mathchardef\measuredangle="0\msx@5D
\mathchardef\sphericalangle="0\msx@5E
\mathchardef\varpropto="3\msx@5F
\mathchardef\smallsmile="3\msx@60
\mathchardef\smallfrown="3\msx@61
\mathchardef\Subset="3\msx@62
\mathchardef\Supset="3\msx@63
\mathchardef\Cup="2\msx@64

\mathchardef\Cap="2\msx@65

\mathchardef\curlywedge="2\msx@66
\mathchardef\curlyvee="2\msx@67
\mathchardef\leftthreetimes="2\msx@68
\mathchardef\rightthreetimes="2\msx@69
\mathchardef\subseteqq="3\msx@6A
\mathchardef\supseteqq="3\msx@6B
\mathchardef\bumpeq="3\msx@6C
\mathchardef\Bumpeq="3\msx@6D
\mathchardef\lll="3\msx@6E

\mathchardef\ggg="3\msx@6F

\mathchardef\circledS="0\msx@73
\mathchardef\pitchfork="3\msx@74
\mathchardef\dotplus="2\msx@75
\mathchardef\backsim="3\msx@76
\mathchardef\backsimeq="3\msx@77
\mathchardef\complement="0\msx@7B
\mathchardef\intercal="2\msx@7C
\mathchardef\circledcirc="2\msx@7D
\mathchardef\circledast="2\msx@7E
\mathchardef\circleddash="2\msx@7F
\def\ulcorner{\delimiter"4\msx@70\msx@70 }
\def\urcorner{\delimiter"5\msx@71\msx@71 }
\def\llcorner{\delimiter"4\msx@78\msx@78 }
\def\lrcorner{\delimiter"5\msx@79\msx@79 }
\def\yen{\mathhexbox\msx@55 }
\def\checkmark{\mathhexbox\msx@58 }
\def\circledR{\mathhexbox\msx@72 }
\def\maltese{\mathhexbox\msx@7A }
\mathchardef\lvertneqq="3\msy@00
\mathchardef\gvertneqq="3\msy@01
\mathchardef\nleq="3\msy@02
\mathchardef\ngeq="3\msy@03
\mathchardef\nless="3\msy@04
\mathchardef\ngtr="3\msy@05
\mathchardef\nprec="3\msy@06
\mathchardef\nsucc="3\msy@07
\mathchardef\lneqq="3\msy@08
\mathchardef\gneqq="3\msy@09
\mathchardef\nleqslant="3\msy@0A
\mathchardef\ngeqslant="3\msy@0B
\mathchardef\lneq="3\msy@0C
\mathchardef\gneq="3\msy@0D
\mathchardef\npreceq="3\msy@0E
\mathchardef\nsucceq="3\msy@0F
\mathchardef\precnsim="3\msy@10
\mathchardef\succnsim="3\msy@11
\mathchardef\lnsim="3\msy@12
\mathchardef\gnsim="3\msy@13
\mathchardef\nleqq="3\msy@14
\mathchardef\ngeqq="3\msy@15
\mathchardef\precneqq="3\msy@16
\mathchardef\succneqq="3\msy@17
\mathchardef\precnapprox="3\msy@18
\mathchardef\succnapprox="3\msy@19
\mathchardef\lnapprox="3\msy@1A
\mathchardef\gnapprox="3\msy@1B
\mathchardef\nsim="3\msy@1C
\mathchardef\ncong="3\msy@1D

\mathchardef\varsubsetneq="3\msy@20
\mathchardef\varsupsetneq="3\msy@21
\mathchardef\nsubseteqq="3\msy@22
\mathchardef\nsupseteqq="3\msy@23
\mathchardef\subsetneqq="3\msy@24
\mathchardef\supsetneqq="3\msy@25
\mathchardef\varsubsetneqq="3\msy@26
\mathchardef\varsupsetneqq="3\msy@27
\mathchardef\subsetneq="3\msy@28
\mathchardef\supsetneq="3\msy@29
\mathchardef\nsubseteq="3\msy@2A
\mathchardef\nsupseteq="3\msy@2B
\mathchardef\nparallel="3\msy@2C
\mathchardef\nmid="3\msy@2D
\mathchardef\nshortmid="3\msy@2E
\mathchardef\nshortparallel="3\msy@2F
\mathchardef\nvdash="3\msy@30
\mathchardef\nVdash="3\msy@31
\mathchardef\nvDash="3\msy@32
\mathchardef\nVDash="3\msy@33
\mathchardef\ntrianglerighteq="3\msy@34
\mathchardef\ntrianglelefteq="3\msy@35
\mathchardef\ntriangleleft="3\msy@36
\mathchardef\ntriangleright="3\msy@37
\mathchardef\nleftarrow="3\msy@38
\mathchardef\nrightarrow="3\msy@39
\mathchardef\nLeftarrow="3\msy@3A
\mathchardef\nRightarrow="3\msy@3B
\mathchardef\nLeftrightarrow="3\msy@3C
\mathchardef\nleftrightarrow="3\msy@
\mathchardef\divideontimes="2\msy@3E
\mathchardef\varnothing="0\msy@3F
\mathchardef\nexists="0\msy@40
\mathchardef\mho="0\msy@66
\mathchardef\eth="0\msy@67
\mathchardef\eqsim="3\msy@68
\mathchardef\beth="0\msy@69
\mathchardef\gimel="0\msy@6A
\mathchardef\daleth="0\msy@6B
\mathchardef\lessdot="3\msy@6C
\mathchardef\gtrdot="3\msy@6D
\mathchardef\ltimes="2\msy@6E
\mathchardef\rtimes="2\msy@6F
\mathchardef\shortmid="3\msy@70
\mathchardef\shortparallel="3\msy@71
\mathchardef\smallsetminus="2\msy@72
\mathchardef\thicksim="3\msy@73
\mathchardef\thickapprox="3\msy@74
\mathchardef\approxeq="3\msy@75
\mathchardef\succapprox="3\msy@76
\mathchardef\precapprox="3\msy@77
\mathchardef\curvearrowleft="3\msy@78
\mathchardef\curvearrowright="3\msy@79
\mathchardef\digamma="0\msy@7A
\mathchardef\varkappa="0\msy@7B
\mathchardef\hslash="0\msy@7D
\mathchardef\hbar="0\msy@7E
\mathchardef\backepsilon="3\msy@7F
\def\Bbb{\ifmmode\let\next\Bbb@\else
 \def\next{\errmessage{Use \string\Bbb\space only in math
 mode}}\fi\next}
\def\Bbb@#1{{\Bbb@@{#1}}}
\def\Bbb@@#1{\fam\msyfam#1}

\catcode`\@=12

\newtheorem{prop}{Proposition}
\newtheorem{lem}[prop]{Lemma}

\newtheorem{cor}[prop]{Corollary}
\newtheorem{them}[prop]{Theorem}

\newtheorem{defnp}[prop]{Definition}
\newtheorem{numconp}[prop]{Construction}
\newtheorem{numrmkp}[prop]{Remark}
\newtheorem{numexp}[prop]{Example}
\newtheorem{numrmksp}[prop]{Remarks}
\newtheorem{conp}[prop]{}

\newtheorem{warningp}{Warning}
\newtheorem{notep}{Note}
\newtheorem{claimp}{Claim}
\newtheorem{examplep}{Example}
\newtheorem{examplesp}{Examples}
\newtheorem{rmkp}{Remark}
\newtheorem{rmksp}{Remarks}

\newenvironment{pf}{\begin{trivlist}\item[]{\sc Proof.}}%
            {\nolinebreak $\Box$ \end{trivlist}}

\newcommand{\noprint}[1]{}

\def\ep{\epsilon}
\def\de{\delta}
\def\nbhd{neighbourhood}

\def\stt{such that}
\def\>>>{\rightarrow }
\def\stt{such that }
\def\al{\alpha}
\def\ep{\epsilon}
\def\x{\times}

\renewcommand{\tilde}{\widetilde}

\newcommand{\GG}{{\frak G}}

\newcommand{\rr}{{\Bbb R}}

\newcommand{\sS}{{\cal S}}

\newcommand{\uU}{{\cal U}}

\newcommand{\vV}{{\cal V}}

\newcommand{\ol}{\overline}

\newcommand{\ldiag}[1]%
       {\makebox[0cm]{${\scriptstyle#1}\downarrow\phantom{\scriptstyle#1}$}}
\newcommand{\ldiagup}[1]%
       {\makebox[0cm]{${\scriptstyle#1}\uparrow\phantom{\scriptstyle#1}$}}
\newcommand{\rdiag}[1]%
       {\makebox[0cm]{$\phantom{\scriptstyle#1}\downarrow{\scriptstyle#1}$}}
\newcommand{\sediagr}[1]%
       {\makebox[0cm]{$\phantom{\scriptstyle#1}\searrow{\scriptstyle#1}$}}
\newcommand{\nediagr}[1]%
       {\makebox[0cm]{$\phantom{\scriptstyle#1}\nearrow{\scriptstyle#1}$}}
\newcommand{\rdiagup}[1]%
       {\makebox[0cm]{$\phantom{\scriptstyle#1}\uparrow{\scriptstyle#1}$}}
\newcommand{\swdiag}[1]%
       {\makebox[0cm]{$\phantom{\scriptstyle#1}\swarrow{\scriptstyle#1}$}}
\newcommand{\sediag}[1]%
       {\makebox[0cm]{${\scriptstyle#1}\searrow\phantom{\scriptstyle#1}$}}
\newcommand{\nediag}[1]%
       {\makebox[0cm]{${\scriptstyle#1}\nearrow\phantom{\scriptstyle#1}$}}

\newcommand{\doublearrowstack}[2]%
{{{{\scriptstyle#1}\atop{\textstyle\longrightarrow}}
\atop{{\textstyle\longrightarrow}\atop{\scriptstyle#2}}}}
\newcommand{\rightleftarrowstack}[2]%
{{{{\scriptstyle#1}\atop{\textstyle\longrightarrow}}
\atop{{\textstyle\longleftarrow}\atop{\scriptstyle#2}}}}
\newcommand{\leftrightarrowstack}[2]%
{{{{\scriptstyle#1}\atop{\textstyle\longleftarrow}}
\atop{{\textstyle\longrightarrow}\atop{\scriptstyle#2}}}}

\newcommand{\overtoparrow}%
{\makebox[0cm]{\beginpicture
\setcoordinatesystem units <.8cm,.4cm> point at 0 0
\setplotarea x from -3 to 3, y from 0 to 1
\setquadratic
\plot -3 0 0 1 3 0 /
\put{\vector(3,-1){0}}[Bl] at 3 0
\endpicture}}

\newcommand{\underbottomarrow}%
{\makebox[0cm]{\beginpicture
\setcoordinatesystem units <.8cm,.4cm> point at 0 0
\setplotarea x from -3 to 3, y from 0 to 1
\setquadratic
\plot -3 1 0 0 3 1 /
\put{\vector(3,1){0}}[Bl] at 3 1
\endpicture}}

\newcommand{\ses}[5]%
{0\longrightarrow#1\stackrel{#2}{ \longrightarrow}#3\stackrel{#4}{
\longrightarrow}#5\longrightarrow0}

\newcommand{\dt}[6]%
{#1\stackrel{#2}{longrightarrow}#3 \stackrel{#4}{\longrightarrow}#5
\stackrel{#6}{\longrightarrow} #1[1]}

\newcommand{\cat}[1]%
{(\mbox{\rm #1})}

\begin{document} \sloppy
\date{October  1999}
\maketitle

\begin{abstract}
Let $X$ be a $G$-space such that
the orbit space
$X/G$ is metrizable. Suppose a family of
slices is given at each point of $X.$ We study a
construction which
associates, under some conditions on the family of slices,
 with any metric on
$X/G$ an invariant metric on $X.$
 We show also that a family
of slices with the required properties
exists for any action of a countable group
on a locally compact and locally connected metric space.

\noindent Keywords: G-space, invariant metric, slice

\noindent AMS classification: 54H15, 57S30

\end{abstract}

\newcommand{\ms}{\ol{M}(W,\tau)}


\section{Introduction}
Let $G$ be a topological group acting continuously
on a topological space $X.$
Existence of a G-invariant metric
on $X$ compatible with the topology of $X$
is a distinct property of the action and provides an important
technical tool.
The obvious necessary condition is metrizability of $X,$
and we assume this in the sequel.
The standard case when an invariant metric
is given by a simple formula, is the case of a compact group  $G$
endowed with a finite measure $\nu$ invariant under (say, right)
translations.
Then from any metric  $d$ on $X$ we obtain an invariant
metric by averaging:
$$ \tilde d(x,y) \ = \ \int_G d(gx,gy)\nu (dg).$$
The most important example is the Haar measure on a compact Lie
group.

For  a topological group $G$
 acting on itself by (say, left) translations, the
 theorem of Birkhoff \cite{B} and Kakutani \cite{Ka} (cf. \cite{MZ})
 gives a sufficient condition
for the existence of an invariant metric, namely
the existence of a countable base of open sets at the unit.
The general case is less known.
Koszul \cite{Ko} and Palais \cite{P} showed that for
 proper (in a sense)  actions of Lie groups
on a separable and metrizable space there
exists an invariant metric.


A similar existence problem can be stated for  riemannian
metrics and actions on smooth manifolds.
Sufficient conditions  for the existence of an invariant
riemannian metric were given by
Koszul
\cite{Ko}, Palais \cite{P} and  Aleeksiejevski \cite{A}.

In the present paper we follow ideas coming from
 differential geometry. In order to explain this in a
 comprehensible way we discuss
first in some detail
the case of totally discontinuous actions (in
particular, the group $G$ is discrete).
This is equivalent to assuming that the
projection $p: X \>>> X/G,$ where $X/G$ is
the space of orbits
 with the quotient topology, is a covering map. Assuming
$X,$ $X/G$ and the projection to be smooth, one can lift
any riemannian metric from $X/G$ to $X,$
since the map $p$ is a local diffeomorphism. 
To measure the distance between points $x,y \in X,$
we have to measure the length of geodesics joining
$x$ with $y.$ This can be done in $X/G,$ since
the projection $p$ becomes a {\bf local isometry}, in particular,
geodesics are preserved  {\bf locally}. The resulting metric is
invariant by construction, since it is defined in terms
of the quotient space.

One can  follow this  observation to construct
an invariant metric for a general covering
space with a metrizable base.
So for any metric $d$ on $X/G$ and a given family $\uU$ of
open subsets of $X$ define

  \begin{eqnarray}\label{for11} \rho (x,y) \ = \ inf \sum_{i=1}^{k-1}
d(px_i, px_{i+1}),\end{eqnarray}

\noindent where $x_1,...,x_k$ is an {\em allowable  sequence}
of points of $X,$
which means that
$x_1 = x,$ $x_k = y$  and any two adjacent points lie in
one  set belonging to  the family $\uU.$
The infimum is taken over all allowable sequences.
Whatever family $\uU$ we are given, we will call its elements
{\it small sets}.
\vskip.5cm
The minimal requirement for $\uU$ is that a small set $U$ should
be elementary, i.e. it  is mapped homeomorphically  by  $p$ onto
its image (which we also call an elementary set).
We also want the family to be $G-$invariant,
to ensure $G-$invariance of $\rho .$
The definition of $\rho$ is the analogue of lifting geodesics
piece by piece when they are divided so that each segment lies
in  an elementary  subset.
By the definition of a covering map, the quotient space has
a base of open elementary sets.
However, if we allow all elementary sets in (\ref{for11}),
then $\rho$ is only a pseudometric. One can see
this easily even in the simplest example of the action
of integers on the real line, which gives a covering
of the circle by $\rr$.
For  locally connected spaces
it is not difficult to find a remedy.
We consider  all balls which are elementary and
balls of radius 4 times larger are still elementary.
The detailed discussion
of that case together with a generalization (to coverings
in the category of spaces with group actions)
is given in Section 3.

\vskip.5cm
We want to explain when  the formula (\ref{for11}) works, hence
in fact we investigate the  problem of lifting a given metric on $X/G$
to $X$ instead of mere existence of an invariant metric.
The assumption of metrizability of $X/G$ is
therefore natural. Moreover, it is easy to observe that
any invariant metric comes from $X/G$ if the latter is
metrizable.

One can observe the strength of the assumption considering
those actions on the real line $\rr$ which have metrizable
orbit spaces.
There are only few possibilities for the latter: the line itself,
the halfline, the closed interval,
the circle and the point. According to this, we may classify
groups which act
effectively on the real line with a metrizable quotient. Any such
action is equivalent to one of  the following:

      i) the (trivial) action of the trivial group;

ii) the involution $x \>>> -x;$

 iii) the action of integers by $x \rightarrow {x+n}$;

iv) the action of the subgroup of isometries of
$\Bbb R$ generated by two elements:
$x \rightarrow {x+1}$ and $x \rightarrow {-x}$ ;

v) a transitive action.
\vskip.5cm

This implies that for any effective action
on $\Bbb R$ with a metrizable quotient  there exists an
invariant  metric except for those examples in case v)
when the group
is in a sense  larger than the orbit
(e.g. the whole affine group).
This shows that the group must be related to the topology
of the orbit, if we want to find an invariant metric .
One consequence of this is that the action must be proper
in a sense (there are  various  meanings of this
word in literature).
We use the name perfect action
(as considered by Koszul and Alekseevski) for
the type of properness we use.
It is not difficult to give examples
of non-proper actions with a metrizable quotient.
We observe
that any perfect action of a countable group
on a locally compact  space has a metrizable
orbits space. This and other topological preliminaries are given
in Section 2.
Then, as a consequence of our construction and some classical
theorems of P.A Smith and Mongomery - Zippin, we show that any
discrete action (cf. Definition \ref{defdis})  on
a manifold is perfect.
\vskip.5cm

The main aim of the present paper is to
extend the construction of the lifting beyond totally
discontinuous actions.

In the presence of nontrivial isotropy (while still assuming
the orbits to be discrete) we
have  to change  formula (\ref{for11})
slightly.
This is because  $p$ can not be a local isometry in general,
hence in an allowable sequence we have to impose a stronger
condition.
The condition corresponds  to the radial distance preserving,
which means that in any small set there is a base point (origin)
\stt the map $p$ preserves the distance between the base point
and any other point of the set.
One can see the problems we encounter even in the simple
example  of the finite cyclic group acting on the
plane $\Bbb R^2$ by rotations. One can see how the properties
of $\rho$ depend on various choices of the class of
small sets.

The next level is the case when the  orbits (hence the group)
are not discrete. Like in the case of locally trivial bundles,
we have to introduce into the formula the metric along
the orbit. Such a family of metrics in orbits
we call an {\it orbital distance}
when some subcontinuity conditions are satisfied (cf. Section 4).
The small sets have to be more refined, in order to remain
related to the topology of the orbit space. The appropriate
requirement is that small sets should be slices
(cf. Definition \ref{defsl}).

The usefulness of slices is clear since the Palais paper \cite{P}
and in fact this notion plays the role analogous to horizontal
distributions in locally trivial smooth fibrations. From
our point of view
this is crucial, as the lift of a geodesic from the base to the
total space of a fibration requires the use of a horizontal
distribution (or a connection).

Having now a family of small  slices $\{ S_x : x \in X\}$ and
an orbital distance $d_O,$ we can again define the lift
of a metric $d$ on $X/G.$

\vskip.5cm

\begin{defnp}\label{defro} We say that $x_1,...,x_k$ is an
{\it allowable} sequence of points in $X$ if for any two
consecutive points one of them belongs to a small slice  at the
other, or they are in the same orbit.

\end{defnp}

Let $\rho$ be the function

\begin{eqnarray}\label{for13}
\rho (x,y) \ = \ inf \sum_{i=1}^{k-1} (d(px_i, px_{i+1})
+ d_O(x_i,x_{i+1})),\end{eqnarray}

\noindent where the infimum is taken over all allowable sequences
\stt
$x_1 = x,$ $x_k = y.$

\vskip.5cm
Our main result consists of two parts. The first is
Theorem \ref{th52} which describes properties of the
orbital metric and
of the family of small sets sufficient to ensure
that the formula (1.2)
gives an (invariant) metric. The properties are independent
of the given metric $d,$ hence once we have such a family,
any metric can be lifted.

Then we give an existence theorem. Theorem \ref{thA}  shows that
for any locally compact and locally connected metrizable
space $X$ and any discrete action of a countable group there
exists a family satisfying the hypothesis of Theorem \ref{th52}.

So the class of spaces for which the construction works is
rather large, as
it contains all manifolds and all locally finite polyhedra.

\section{Topological  preliminaries}

{\bf Notation}

Given a metric $\rho$, by $K_{\rho}(x,\ep)$ we denote
the ball of radius
$\ep$  centered at $x$ with respect to $\rho.$

In the sequel, all actions will be assumed continuous.
 For an action of a group $G$ on $X$ we denote
by $X/G$ the orbit space and by $p: X \>>> X/G$
the quotient map. By $G_x$ we denote the
stabilizer (isotropy subgroup)
of $x;$
$G_x =\{ g\in G: gx=x\}.$

$GA = \{ ga : g \in G, a \in A\} ,$ hence $Gx = G\{ x\}$ denotes
the orbit of $x.$

The map  $e_x : G \>>> Gx ; g \>>> gx$ will be
called the {\it evaluation map.}

\begin{defnp}\label{defsl} By a slice at $x \in X$ we
mean a subset $S_x$

\stt\

1.  $x \in S_x;$

2. $gS_x \cap S_x \neq \emptyset$ implies $gx = x;$

3. $gS_x = S_x$ for $g \in G_x,$ hence $G_x$ acts on $S_x;$

4. the map $(G \x S_x)/G_x \>>> X$ given by the evaluation
is a homeomorphic embedding
onto an open neighbourhood of the orbit $Gx.$
 \end{defnp}

Note that our definition is stronger than the usual one, since
it includes
condition 4, which implies that the evaluation map $G \>>> Gx$
is open for any $x.$ Under this assumption
the topology of $G$ is strictly related to the topology
of orbits. In particular, this excludes the examples given
by changing the topology of a group acting on $X,$ for instance
from a connected group to a discrete group.
Note the following properties of slices:

1) the image of any slice is open in $X/G;$

2) for any open set $U \subset X/G,$ the
intersection $p^{-1}U \cap S_x$
is open in $S_x$ and, if $px \in U,$  it is another slice at $x.$

For an action of a discrete group, the slice is an open subset.
This leads to the following definition.

\begin{defnp}\label{defdis} The action of a
(discrete) group $G$ on a topological space
$X$ is called discrete if for any orbit
$Gx$ there exists an open neighbourhood $U$ of
$x$ \stt for  all $h\not\in G_x$ we have
 $hU \cap U = \emptyset$
 and every  point $x\in X$
has a base of open $G_x$-invariant neighbourhoods.
\end{defnp}

Any orbit of a discrete action is
discrete in
$X$.
If the action is free and discrete, then it is  totally
discontinuous.
\medskip

The following observation is well known (cf. [P], p. 319).

\begin{prop}\label{pal}
Let $ \Gamma$ be a group acting by isometries on a
metric space $(X,\rho)$ such that
each orbit $ \Gamma{x}$ is closed in $X.$
Then the formula
$$ d(px,py)  =
inf \{ \rho (a,b) :a \in  \Gamma{x},
b \in  \Gamma{y} \}$$
defines a metric on the orbit space.
\end{prop}

In the sequel we will need the following fact.
If $X/G$ is  $T_1,$ then  any orbit is a closed subset of $X.$
Thus any accumulation point of $Gx$ belongs to that orbit.
Therefore either there are no accumulation points or,
by homogeneity,
every point of
$Gx$ is its accumulation point. When $X$ is locally compact,
in the latter case  we would get an uncountable orbit.
This yields the following lemma.

\begin{lem}  Let $G$ be a countable group acting on
locally compact Hausdorff space $X$
such that the quotient space is a $T_1$-space. Then any orbit
is discrete in $X.$
\end{lem}

When we want to apply our construction,  the following
notion arises naturally. It is stronger  than   the usual notion of
proper actions.

\begin{defnp}\label{defprop} We call an effective action
of a group $G$ on a
Hausdorff space $X$ perfect
if for any compact sets $U,V \subset X$
the set $ \{ g \in G : gU \cap
V \neq \emptyset \}$
is compact.
\end{defnp}

Any perfect action has compact stabilizers (hence
finite, if the group is countable). If a  countable group acts
perfectly on a locally compact space $X$, then any orbit
is discrete in $X.$

Now we will show that the quotient space of any perfect action
of a countable discrete group on
a locally compact metrizable space is metrizable.
Note, however, that perfectness is not a necessary
condition.  The following example
shows a nonperfect action with a decent quotient.

Consider the  subset
$Y  = \{(x,0) : x \in \Bbb R \} \cup \bigcup_n
Y_n$  of the plane, where $Y_n$ is the
sum of two closed intervals with a common initial point
$(n,0)$, for every $n \in \Bbb Z.$
The  group acting on $Y$   is generated by homeomorphisms
 $g_n$ which exchange the two intervals attached at $(n,0)$
and fix the rest of $Y$.

\begin{them}\label{th31}     Let $G$ be a countable
group acting on a locally compact
metrizable  space $X$.
If  the action is perfect, then the group is discrete and
the quotient space $X/G$ is
metrizable. \end{them}

\begin{pf}
First note that by  the remark above,
each orbit
$Gx$ is closed, so $X/G$ is a $T_1$ space.
Note also that under the assumptions of the theorem, for
any compacts
$K_1,K_2 \subset X$ the
set $\{g \in G : gK_1 \cap K_2 \neq \emptyset \}$ is finite.

     Now let $A$ be a closed subset of $X/G$ and $Gx_0$ any
orbit which does not belong to $A$.
Then there exists $\delta   > 0$ such that
$\overline{K(z,\delta)}$ is compact and
$\overline{K(z, \delta)} \cap p^{-1}(A) = \emptyset$,
where $z \in Gx_0$.

For any $y \in p^{-1}(A)$
one can find $K(y, \de_y)$ which is relatively  compact
and \stt  $Gx_0
\cap \overline{K(y,\de_y)}  =  \emptyset .$
 Then
$$ B  = \{ g \in G : g(\overline{K(z,\delta)}) \cap
\overline{K(y,\de_y)} \neq \emptyset \} $$
is finite, so there exists $\epsilon_y \leq \de_y$ such that
$K(y,\epsilon_y)$ is disjoint with
$\bigcup_{g \in G} g\overline{K(z, \delta)}.$
Since $(\bigcup_{g \in G} g(K(z,\delta))) \cap (\bigcup_{g \in G}
g(K(y,\epsilon_y)) =  \emptyset ,$
we see that the space $X/G$ is regular.

Now the following result of A.H. Stone (cf. \cite{E}, 4.5.17)
completes the proof.
\end{pf}

{\begin{them}  If $f:X \rightarrow Y$ is a continuous and
open mapping from a  locally
separable metric space $X$ onto a regular space $Y$ and the set
$f^{-1}(y)$ is separable for
every $y \in Y$, then $Y$ is a metrizable space.
\end{them}}


 \section{Covering spaces}

In this section we consider totally discontinuous
actions of discrete groups, so that the projections
$p : \tilde X \>>> X = {\tilde X}/G$ are coverings.
Our aim is twofold. First
we want to explain our construction in a relatively simple case.
Secondly we show the existence of a lifted metric  for
coverings in the category of spaces
with group actions.

\begin{defnp}\label{defca} Consider a group epimorphism
$\phi:{\tilde G} \rightarrow G.$ Suppose that a group
$G$ acts on $X$ and $\tilde G$ acts on $\tilde X.$
By a
covering in the category of spaces
with group actions we mean  a covering map
$p: \tilde X \>>> X$ \stt

$$p({\tilde g}(\tilde x))  =   \phi{({\tilde g})}(p{\tilde x})$$

\noindent for any $\tilde x \in \tilde X$ and
$\tilde g \in \tilde G$.
\end{defnp}
In the following we assume that $\tilde X$ is connected. If not,
one can apply our construction to every connected
component of $\tilde X.$

\begin{defnp}\label{defele} Let $p:{\tilde X}\rightarrow X$
be a covering map.
 Any open set $U \subseteq X$ \stt its counterimage
$p^{-1}(U)$ is homeomorphic to a disjoint sum
$\bigcup _{j\in J}U_j$ and $p$ restricted
to $U_j$ is a homeomorphism onto $U,$ is called elementary.
Any of its homeomorphic copies $U_j$ will be
also called elementary.
\end{defnp}

\begin{them} For any covering map $p:{\tilde X}\rightarrow X$
of a localy connected space $X$
and for any metric $d$ on $X$ there exists a metric
$\rho$ on
$\tilde X$ such that $p$ is locally isometric.
Moreover, if $X$ is a $G$-space
and $\tilde X$ is a $\tilde G$-space covering of $(X,G)$,
then any $G$-invariant
metric lifts to a $\tilde G$-invariant metric on $\tilde X$.
\end{them}

\begin{pf}
Since $X$ is locally connected, thus each point of $X$
(and thus of $\tilde X$)
has a base of connected elementary
neighbourhoods.

      Let $x \in X$. Take a connected elementary
neighbourhood $U_{x}$
of $x.$ Then there exists a ball
$K(x,\epsilon_x) \subseteq U_x$
and a connected open set $V_x \subseteq
 K(x,\epsilon_{x/4})$ containing $x$. Define
now a family $\vV$ of small sets as the family of all
connected components of  sets $p^{-1}V_x$ for all $x\in X.$
Note that all such  small sets are elementary, since all $V_x$
are.
As in Definition \ref{for13} define now
\vskip.5cm
$$\rho{(a,b)}  =
inf\, \{\sum_{i=1}^{k-1}  d(p\tilde x_i,
p\tilde x_{i+1})\},$$
where the infimum is taken over all
$\vV$-admissible sequences
$\xi = \{ a=\tilde x_1,
\tilde x_2, ... ,\tilde x_{k}=b\}$ in $\tilde X,$ i.e.
such that any two consecutive points of
the sequence  belong to the same small set.

 It is easy to check that $\rho$ is a pseudometric
and $0 < d{(pa,pb)} \leq  \rho{(a,b)}$ when $pa \neq pb$.

      Consider two sets $A,B \in \vV$ such
that $A \cap B \neq \emptyset.$
By the definition of $\vV ,$ we have $pA = V_x, pB = V_y$ for
some $x,y \in X.$
When we assume  that $diam(V_y) \leq diam(V_x)$, then
$diam(V_y) \leq diam(V_x) \leq \epsilon_x/2$ and
$V_x \cap V_y \neq \emptyset,$
so we see that
 $V_x \cup V_y \subseteq K(x,{\epsilon_{x}}) \subseteq U_x.$
Thus $A \cup B$ is
elementary and because of the connectedness
of $V_x$,$V_y$ and $U_x$ it follows
that $A \cup B$ is included in some component of  $p^{-1}U_x$.

    Now take any points ${\tilde x},{\tilde{\tilde x}}$ such that
$p{\tilde x}  =  p{\tilde{\tilde x}}  =  x$.
 We will show that every admissible
sequence joining $\tilde x$ and $\tilde{\tilde x}$ contains
$z \not\in p^{-1}(V_x)$. Assume the contrary.
Then in such an
admissible sequence there exist ${\tilde x_i},{\tilde x_{i+1}}$
which belong
to different components of $p^{-1}(V_x)$, and  belong
to a component of
$p^{-1}(V_z)$ for some $z \in X$. We have now two cases:

      i) $V_x \cup V_z \subseteq S_x$. Then some
      component of $p^{-1}(S_x)$ intersects
two different components of $p^{-1}(V_x)$, which is impossible.

     ii) $V_x \cup V_z \subseteq U_z$. Then each component
     of $p^{-1}(V_x)$ is
contained in precisely one component of
$p^{-1}(U_z)$ ($p:{\tilde X} \rightarrow X$
is a covering map), but at least one component of $p^{-1}(V_z)$
intersects two
different components of $p^{-1}(V_x)$, which is a contradiction.
\vskip.5cm
      Because $V_x$ is open, so
$0 <   \delta \leq \rho{({\tilde x},{\tilde{\tilde x}})}$.
Note also that $p$ restricted to any component of any
admissible set is an
isometry. This yields compatibility of $\rho$ with the
topology of $\tilde X.$


To prove the second part of the theorem take as small sets
the family of the connected components of the sets
$\{ p^{-1}(gV_x) : g \in G , x \in X \}$.
      First note that if $U$ is elementary then $gU$ is also
elementary, because the action of $G$ is covered by the action
of $\tilde G.$ Thus $\vV$ is $\tilde G$-invariant.
Since $G$ acts by isometries, so  $\rho$ defined as before is
$\tilde G$-invariant. By the same argument  as in part
one of the theorem
we  obtain that $\rho$ is a metric on $\tilde X$
compatible
with the topology on $\tilde X$. \end{pf}


\section{Orbital distance}

When a $G$-space is endowed with an invariant metric, then
we have
in particular a continuous family of metrics in orbits.
On the other hand, an important  ingredient of  our construction
is a family of invariant  metrics in orbits. We follow
the construction which works for locally product bundles.
In this case, with help
of a horizontal connection, a metric on the total space can
be obtained
from a metric on the base and a fiberwise family of metrics.

We encounter here two main difficulties. Even when the metric
on an orbit
is induced from a metric on the group $G,$ then there is
no canonical way to do this, because the identification of the
orbit with the quotient $G/G_x$ requires an apriori
choice of the base point $x$ in the orbit. Once an orbit has a
base point $x,$ one may want to use the slice $S_x$
to endow the orbits
in the \nbhd\  $GS_x$ of $Gx$  with base points.
But neighbour orbits in general cut the slice
at many points, so we are forced to assume that
any such base point gives
the same metric (i.e, that  the metric on $G$ is right
invariant with respect to $G_x).$
The other problem is that the family of metrics coming from
$G$ is not continuous in general. It has, however, some
subcontinuity
properties which are sufficient for our purposes.

Let $d_G$ be a left invariant metric on $G.$
Consider first a single orbit with a base point $x$ in it and
 assume that the evaluation map
$e_x : G/G_x \>>> Gx$ given by  $e_x(g) = gx$
is a homeomorphism.
Using this identification, we want to endow the orbit with
an invariant metric. For a subgroup $K \subset G$
define  a left $G$ -invariant pseudometric on $G/K$
by the formula

\begin{eqnarray}
\label{defhom}  d_K(g_1K, g_2K) =
inf \{ d_G(g_1u, g_2v) : u,v \in K\} .
\end{eqnarray}

In particular,
$$ d_K(g_1K, g_2K) \leq d_G(g_1,g_2)$$
for any $g_1, g_2 \in G$

\begin{lem}\label{exor}
Let $H$ be a closed subgroup of the group $G.$
If $d_G$ is left $G$-invariant and right $H-$invariant,
then for any closed subgroup $K\subset H$ we have
\begin{eqnarray}\label{exor1}
d_K(g_1K,g_2K) = inf \{ d_G(g_1, g_2u) : u\in K\},
\end{eqnarray}
thus formula (\ref{defhom}) defines a
$G$-invariant metric on $G/K,$ compatible with the quotient
topology.
The same is true if $K$ is a normal subgroup of $G.$
When we replace  $K$ by $K' = hKh^{-1}, h\in H,$
then $G/K$ and $G/hKh^{-1}$ with the
metrics $d_K, d_{hKh^{-1}}$ are isometric under the
natural
map $$\phi (gK) = gh^{-1}K'.$$
\end{lem}

\begin{pf}
Both  formula (\ref{exor1}) and the isometry follow
directly from the right
invariance.
For instance, we have
$d_K(g_1K, g_2K) = d_K(g_1h^{-1}hK, g_2h^{-1}hK)
= d_{hKh^{-1}}(\phi (g_1K), \phi (g_2K)).$
Since we assumed $K$ closed, thus $d_K $ is a metric.
\end{pf}

Under the identification of the  orbit $Gy$ with $G/G_y$
by the evaluation map, the metric $d_{G_y}$ is a well-defined
metric in the orbit.
The translation of the
orbit by $h \in G_x$ which changes
the point $y$ to $hy$ corresponds to
the above homeomorphism $\phi :G/G_y \>>> G/G_{hy}.$
Given a slice $S_x$ at $x$ and a metric in $G$ which is
left $G-$invariant and right $G_x-$invariant we have a well
defined
metric in any orbit passing through $S_x,$ since we can use
any of the points of intersection $Gy \cap S_x$ as  the
base point
of $Gy$ (the stabilizer $G_x$
corresponds to $H$ and $S_y$ to $K$).
We denote the whole family of metrics by $d_x$
to stress which slice was used to define the metrics.

\begin{cor}\label{sub1}
Let $d_G$ be an invariant metric on $G.$ Given a point   $x$,
assume either that $d_G$ is  right $G_x$-invariant
or that all stabilizers of points in $S_x$ are normal.
Then we get from $d_G$  a
family $d_x$ of $G$-invariant metrics  in all  orbits
intersecting $S_x.$ For any point
$y \in S_x$ and   $g_1, g_2 \in  G$ we have the inequalities
\begin{eqnarray}
\label{sub2}  & d_x(g_1x, g_2x) \leq d_x(g_1y, g_2y) \leq d_G(g_1,g_2).
\end{eqnarray}
If $y \in g_0S_x$ for some $g_0 \in G,$  then
\begin{eqnarray}\label{sub4}
d_x(y, gy)  \leq d_G(g_0, gg_0).
\end{eqnarray}
\end{cor}

\begin{pf}
To prove (\ref{sub2}) we start from the inequality
$d_x(x, gx) \leq d_x(y, gy)$ which holds, by the definition
of $d_x,$ for any $y\in S_x$ and any $g.$
Then we use the left invariance of $d_x.$

Inequality (\ref{sub4}) follows from (\ref{sub2}).
\end{pf}

Thus, assuming that the action has slices, our construction works
locally
for the following classes of actions :

$\bullet$ any action of a discrete group, or more generally
of a group having
a biinvariant metric,

$\bullet$ perfect actions, and in general for actions with compact
 isotropy subgroups,

 $\bullet$ actions whose all stabilizers are normal in $G.$

Now we globalize the procedure by gluing together the metrics
defined locally with the use of a decomposition of unity.
Let a $G$-space $X$ be given  having a slice $S_x$ at each
point such that $S_{gx} = gS_x$ and  suppose  that
$X/G$ is metrizable and endowed with a metric $d.$

For any open set $U\subset X/G,$ the set $p^{-1}U\cap S_x$ is
also a slice at $x,$ and we say that it is a subslice
(with respect to the given family of slices).

Under our assumptions, there exists a family $S_\al,$
each $S_\al$ a subslice
of $S_{x_\al}$ for a point $x_\al \in X,$  \stt
$\overline S_\al \subset S_{x_\al}$ and
the family $pS_\al$ is a locally finite open cover of $X/G.$
The cover admits a decomposition of unity $\{\chi_\al\} .$
Denote by $d_\al$ the local family of metrics in orbits
in $GS_\al $ and define the global family of metrics in orbits
by the formula

\begin{eqnarray}
\label{defo} d_O(x,gx) = \sum_\al \chi _\al(px)\, d_\al (x,gx).
\end{eqnarray}

To describe it formally, denote
$V = \{ (x,y):y=gx, g\in G\} .$
Then $d_O: V \>>> \Bbb R$ is G-invariant and restricted to
each orbit it is a metric. In the sequel we will consider the map
as extended by zero map to the
whole $X\times X$ and we will call it an {\it orbital metric.}

The orbital metric  $d_O$  is
 in each orbit  compatible with the topology
of the orbit (since each metric $d_\al$ is). However, it is not a
continuous function even on $V,$ since orbits with smaller
isotropy
do not converge to  orbits of larger isotropy. For
instance, $d_O$ is not continuous at fixed points.

In order to get an orbital metric with nice properties,
we   assume that the given family of slices
$\{ S_x\}_{x\in X}$
satisfies  the following condition:
\medskip

\noindent (*)\ \ \  for any $x \in X$ and $y\in S_x,$ the
intersection $S_x \cap S_y$ is open in $S_y.$

\medskip

This implies in particular that  any $y \in S_x$ has  an open
neighbourhood
$U$ \stt $S_y \cap U \subset S_x.$

The following conditions are stronger than (*).

(1)  For any $y \in S_x$ we have $S_y \subset S_x.$

(2)  For any $y \in S_x$ and $g \not\in  G_x$ we
have  $S_y \cap S_{gx} = \emptyset .$
\medskip

First we  show that assuming (*) we have some subcontinuity
properties for $d_O.$
Let $S_x(\de ) = S_x \cap p^{-1}K_d(px,\de ).$

\begin{conp}\label{propa} {\bf Property A.} For any $x \in X$ and
any $\ep > 0$ there exists $\de > 0$
\stt $d_O(y, gy) < \ep$ if $y \in S_x(\de )$ and $d_G(g,1)<\de .$
\end{conp}

\begin{pf}
A given point $x$ is in the closure of a
finite number  of sets $GS_\al ,$ say $GS_1,...,GS_r.$
It means that $x \in GS_{x_1}\cap ...\cap GS_{x_r},$ where
$S_{x_i}$ is the slice which by our assumption
contains the closure of $S_i.$ Let $g_i$ be chosen
\stt $g_ix \in S_{x_i}.$ By (*)
there exists $\de > 0$ \stt for
all $y \in S_x(\de )$ we have $g_iy \in S_{x_i}.$
We can also assume that $S_x(\de )$ does not
intersect $S_\al $ except for $S_1,...,S_r.$
Thus
$$\begin{array}{cc}
d_O(y,gy) = & \sum_{i=1}^r \chi _i(px_i)\, d_{x_i}(y, gy) = \\
            & \sum_{i=1}^r \chi _i(px_i)\,
            d_(g_iy, g_igg_i^{-1}g_iy)\leq\\

              & \sum_{i=1}^r \chi _i(px_i)\, d_G(1, g_igg_i^{-1})
\end{array}$$

\noindent by (\ref{sub2}).
Since the conjugation by $g_i$ in $G$ is continuous,
there exists  $\de$ with the required properties.
\end{pf}

Now we will show that the orbital distance is
subcontinuous.

\begin{conp}\label{propb}{\bf Property B.}
For any point $x$ there exists $\de > 0$ \stt
 for any $y \in S_x(\de )$ and any $g_1,g_2 \in G$ we have
\begin{eqnarray}\label{subo2}
d_O(g_1x,g_2x) \leq d_O(g_1y,g_2y).
\end{eqnarray}
\end{conp}

\begin{pf} This is straightforward from the definition of $d_O$
(and from the fact that isotropy is subcontinuous)
once we know (as in A) that the
element $g_i \in G$ which brings $x$ into $S_i$
 does the same with any element of $S_x(\de ).$ But this is true
for $\de $ small enough by (*).
\end{pf}

Similarly,  the assumption that the evaluation map
$G/G_x \>>> Gx$ is a homeomorphism gives the following result.

\begin{conp}\label{propc} {\bf Property C.}
For any point $x\in X$ and any $\de > 0$ there exists
$\ep$ \stt $d_O(x,gx) < \ep$ implies
that there exists $u\in G_x$ satisfying the inequality
$d_G(1,gu) < \de .$
\end{conp}

Given a $G-$space with slices and an invariant
metric $d_G$ in $G,$
by an {\it associated orbital metric} $d_G$ we mean
an orbital metric which
in each orbit  restricts to a
 $G$-invariant metric isometric to $d_K$
for a stabilizer $K$ in the orbit
and has  properties A,B and C.

From the above considerations we know that an   orbital metric
exists
for any $G$-space with slices \stt there exists
a metric on $G$ which is left $G$-invariant and right
$G_x$-invariant
for any stabilizer $G_x$, or
the stabilizers are normal, and
the quotient $X/G$ is metrizable.

\begin{rmkp} The metric in $X/G$ is not necessary in
our construction,
it is enough to assume the quotient is paracompact. More
difficult
it is to show that the condition (*) can be removed when the
quotient is locally compact. We will not go into the
details,
since  we treat here the general case. \end{rmkp}


\section{Main results}

   We will now prove the main theorems of this paper.

Consider an action of a group $G$ on $X$ and suppose that $G$
admits a left invariant metric $d_G.$
Assume  we are given a  family ${\sS} =\{ S_x\} _{x\in X}$
of slices in $X$
\stt $S_{gx} = gS_x$ for any $x, g,$ an associated
orbital metric $d_O$ in $X$
and a metric $d$ in $X/G.$

\begin{defnp}\label{defal}
A sequence   ${x_1,...,x_n}$
is called an ${\sS}$-{\it allowable} (or simply
{\it allowable)} sequence   if
for each $1\leq i \leq {n-1}$ we have $x_i \in S_{x_{i+1}},$  or
$x_{i+1}\in S_{x_i},$ or $x_{i+1} = gx_i.$

For an allowable sequence $\xi = \{ x_1,...,x_n\}$ in $X$ denote
$$\Sigma_{\xi} = \sum_{i=1}^{n-1} d(px_i,px_{i+1}),$$
$$\Phi_{\xi} =  \sum_{i=1}^{n-1} d_O(x_i,x_{i+1}). $$
\end{defnp}

\begin{defnp}\label{defmet}

 $$\rho(x,y) = inf \{ \Sigma_{\xi} + \Phi_{\xi}
\ : \xi  = \{ x_1,...,x_k\} \  {is \ allowable,} \ \ x = x_1, \ y = x_k\} . $$
\end{defnp}

We say that a $G$-space {\it has small slices } if there exists
a slice at every point $x\in X$ and for   any open \nbhd\
of $x$ and any slice $S_x$ there exists a
slice $S_x' \subset S_x\cap U.$

\begin{them}\label{th52}
 Let $G$ be a topological
group with a left invariant metric
$d_G$, acting on $X.$
Suppose  that $X$ has small slices and we  are given  a family
of slices
 $\sS = \{ S_x\} _{x\in X}$ satisfying

(i) $gS_x = S_{gx}$ for any $x \in X$ and $g \in G;$

(ii) for any $y\in S_x,$ \  $S_y \cap S_{gx} \neq \emptyset $
implies $gx = x.$

Suppose also
a metric $d$ in  the orbit space $X/G$ and an orbital distance
$d_O$ are given. Then formula (\ref{defmet}) defines a lift of
the metric $d$ to the
$G-$invariant metric $\rho$ in $X,$  compatible
with the topology of $X.$
\end{them}

\begin{pf}

It is a direct consequence of the definition that $\rho$
is a pseudometric in $X$ and $\rho (x,z) \geq d(px,pz)>0$
if $px \neq pz.$ It remains to prove that $\rho $ distinguishes
points in the same orbit and that it is compatible with the
topology of $X.$
For any $\ep > 0$  we denote $S_x(\ep ) =
S_x \cap p^{-1}K_d(px,\ep )$  and by $B(\ep ) $ the ball in $G$
with the center at the unit element and of radius $\ep .$
Recall that $S_x(\ep )$ is again a slice at $x.$

\begin{lem}\label{lem52}
For any point  $x \in X$ and $\ep > 0$ there exists $\de > 0$ \stt
$$ B(\de )S_x(\de ) \subset K_\rho (x,\ep ).$$
\end{lem}

\begin{lem}\label{lem53}
For any point  $x \in X$ and $\de >0$ there exists a positive
number $\ep$
\stt  $K_d(px,\epsilon ) \subset pS_x$  and
$$ K_\rho (x,\ep ) \subset
B(\de )S_x(\ep).$$
\end{lem}

Proof of Lemma \ref{lem52}.
Consider an allowable sequence $\xi = \{ x, y, gy\} ,$
where $y \in S_x(\de )$ and $g \in B(\de ).$
By Property A of orbital metrics
for  $\de$ small enough, $\rho (x, gy) \leq d_O(y, gy) +
d(px, py) < \ep .$
\medskip

Proof of Lemma \ref{lem53}.
Since the projection $p$ is open, there exists $\ep$ \stt
$K(px,\ep) \subset pS_x.$
Let $z \in K_{\rho}(x,\ep),$
$\xi = \{x=x_1,x_2,...,x_k=z\}$ be an allowable sequence
 and $\Sigma_{\xi} + \Phi_{\xi} < \ep .$
 If an element $x_i$ of the sequence $\xi$ does not
belong to $ GS_x(\ep ), $ then $\Sigma_
\xi \geq d(px, px_i) \geq \ep .$
Thus any point $x_i$ in the sequence $\xi$
can be written as $g_ix_i'$ for some
$g_i \in G, x_i' \in S_x(\ep ).$
We can choose $x_1',...,x_k'$ and  $g_1,...,g_k$ \stt
$x_i' = x_{i+1}'$ when  $x_{i+1} = gx_i$ and $g_i = g_{i+1}$
if $x_{i+1} \in S_{x_i}.$
In the last case we use the assumption (ii), which is equivalent
to the inclusion $S_y\cap GS_x \subset S_x$ for any $y \in S_x.$
It gives that for $x_{i+1} = g_{i+1}x_{i+1}'$ we have
$g_i^{-1}g_{i+1}\in G_x,$ so we have
$x_{i+1} = g_ig_i^{-1}g_{i+1}x_{i+1}'.$

By Property B of the orbital metric, if $\ep$ is small
enough then
$$d_O(x_i,x_{i+1}) = d_O(g_ix_i',g_{i+1}x_i')
\geq d_O(g_ix,g_{i+1}x)$$
in the case $x_{i+1} \in Gx_i,$ and we have
zeros on both sides  otherwise.
It follows that
$$d_O(x,g_kx) \leq \sum_{i=1}^{k-1}d_O(g_ix,g_{i+1}x) \leq \Phi_\xi
< \ep,$$
and by Property C  for $\ep$ small enough
there exists $u \in G_x$
\stt $d_G(1,g_ku) < \de $.
So we have $z = g_kx_k' = (g_ku)(u^{-1}x_k').$
This completes the proof of Lemma \ref{lem53}.

The last lemma implies that $\rho$ distinguishes points
in orbits:
if $\rho (x,gx) = 0,$ then there is a sequence $g_n \in G$
convergent to the unit of $G$ \stt $gx = g_nx.$
Thus $gx = x.$
Since the sets
 $B(\de ) \x S_x(\ep )$ generate the topology of $X,$
Lemmas \ref{lem52} and \ref{lem53} show that $\rho$ is
compatible with the topology
of $X.$
\end{pf}

\begin{rmkp} Assumption (ii) can be omitted  if the quotient
space is locally compact (cf. Section 4). If $X$ is locally
compact, then the assumption of existence of small slices
is superfluous.

\end{rmkp}


When one wants to find  a class of actions for which
the assumptions
of (\ref{th52}) are satisfied, the first try would be  discrete
actions. In that case the orbital metric is trivially given
by the (biinvariant) discrete metric on $G,$ and we have
an easy description of slices. A slice at $x$ is an open \nbhd\
$U$ of $x$ which is $G_x$-invariant and disjoint with
any of its translations by any $g\not\in G_x.$
We prove that for any locally compact and locally connected
$G$-space one can find the family of slices for actions
with a Hausdorff orbit space.
Note that the condition (*) is easy to obtain once we
have slices, because for any slice $S_x,$ any  $ y \in S_x$
and a slice $S_y,$ the set $S_x \cap S_y$ is a slice at $y.$

\begin{them}\label{thA}
Let $X$ be a metrizable, locally compact and
locally connected space and let
$G$ be a countable group acting on $X$ such that the
quotient space $X/G$
is Hausdorff. Then  each orbit $Gx$ is discrete in $X$
and every point $x \in X$ has a base of  connected,
$G_x$-invariant open neighbourhoods.

Moreover, there exists a family $U_x : x\in X$ of
$G_x$-invariant, path
connected open sets satisfying

i) $U_{gx}  =  gU_x$,

ii) if $U_x \cap U_{gx} \neq \emptyset$ then
$g \in G_x.$
\end{them}

\begin{pf} Take any $x \in X$ and an open set $U$ containing $x$.
There exists $V \subseteq U$
such that $x \in V$, $\overline V$ is compact and
$\overline V \cap Gx  = \{x\}$.
We claim
that there exists an open set $W$ \stt
$x \in W  \subseteq \bigcap_{g \in G_x} gV.$
\vskip.5cm
 Assume it is not true. Then for any natural  $n$ there exist a
connected neighbourhood $V_n \subseteq V$ of $x$ such that
$x \in V_n \subseteq K(x,(1/n))$
and  $g_n \in G_x$ such that
$(X - \overline{g_nV}) \cap V_n \neq \emptyset$.

\vskip.5cm
We want to show that $\partial (g_nV) \cap V_n \neq \emptyset.$
If this holds,
then $V_n$ can be decomposed as the sum of the sets
$V_n \cap g_nV$ and $V_n \cap (X - \overline{g_nV})$ which are
disjoint, open and non-empty, so $V_n$ would be disconnected.

Pick any $d_n \in \partial (g_nV) \cap V_n.$ Notice that
$\rho(d_n, x) < \frac1n ,$ where $\rho $ denotes
an auxiliary metric on $X.$ It follows that $w_n = g_n^{-1}d_n
\in \partial V.$ The last set is compact, so up to passing
to a subsequence we may assume
$w_n \>>> w' \in \partial V.$ Since
$\overline V \cap Gx = \{ x\} $ and
$X/G$ is Hausdorff, we get a contradiction.

  An open set $W$ contains a connected  neighbourhood $W'$
of $x$. Now
$ W''  =  \bigcup_{g \in G_x} g{W'}$ is an open, path connected and
$G_x$-invariant subset of  $V.$  Since $V$ can be chosen
arbitrarily small, the proof of the first part  is complete.

    Now the last part of the theorem.

     Notice first that if a neighbourhood $U_x$ of $x$
is $G_x$-invariant
then $U_{gx}  =  gU_x $ is $G_{gx}$-invariant.

     Assume, on the contrary, that for every $G_x$-invariant
     open set $S_x $
containing $x$
there exists $g \not\in G_x$ such that $gU_x
\cap U_x \neq \emptyset$.
\vskip.5cm
     Choose then a $G_x$-invariant neighbourhood $U$ of $x$
such that $\overline U
\cap Gx  = \{x\}$.
For every ball $K(x,(1/n)) \subseteq U$ there exist open
$V_n \subseteq K(x,(1/n))$
which are connected and $g_n \not\in G_x$ such that
$g_n{U} \cap V_n \neq \emptyset$.
The same argument as before
yields a contradiction.
\end{pf}

\begin{them}\label{thB}
 Let $X$ be a locally compact and locally path connected  metric
space and let $G$ be a countable group acting on $X$ such that the
quotient space
$X/G$ is a metric space with a metric $d$. Then there exists
an open covering
$\{S_x\}_{x \in X}$ of $X$ such that

A) $U_{gx}  =  gU_x$,

     B) if $U_x \cap U_{gx} \neq \emptyset$ then $g \in G_x$,

     C) for any $y \in U_x$, if $U_y \cap U_{gx} \neq \emptyset$
then $g \in G_x$
     for any $x \in X$ and $g \in G$.
\end{them}

\begin{pf}
     We already know that a covering satisfying A) and B) exists.
Further, $p:X \rightarrow X/G$ is open, so $pS_x$ is open. Hence
there exists
$K(px,{ \epsilon_x}) \subseteq pS_x$. Now $S_x \cap
p^{-1}(K(px,({ \epsilon_x}/8)))$
is open and $x$ belongs in it, so there exists a
$G_x$-invariant, path connected and
open set $V_x$ contained in the set mentioned above. We will show that
the covering
$\{V_x\}_{x \in X}$ satisfies C), because it obviously
satisfies A) and B).

 Note first  that  for $y \in V_x$ the conditions $V_x
\cup V_{gx} \subseteq V_y$
for some $x,y \in X$ and $g \not\in G_x$ can not be  satisfied
simultaneously. To see this take
 $y \in S_x;$ then
$S_y \subseteq G_x$ because of B); but
also $x \in V_x \subseteq V_y$, so
$G_x  =  S_y$. Now $V_{gx} \subseteq V_y$ for some $g \not\in
G_x$ gives a contradiction by the argument
described above  for  covering spaces.

Assume, on the contrary, that $V_y \cap V_{gx} \neq
\emptyset$ and $y \in V_x$
for some $x,y \in X$ and $g \not\in G_x$.
     We have two cases:

     i) $diam(pV_y) \leq diam(pV_x)$.
     Then $pV_y \subseteq pS_x$ and for
each $g \in G$ there exists $h \in G$ such that
$V_{gy} \subseteq U_{hx}$.

Then it is
a clear contradiction because of A), B) and the connectedness
of all sets considered.

ii) $diam(pV_x) \leq diam(pV_y)$. Then $pV_x \subseteq pU_y$ and for
every $g \in G$ there exists $h \in G$ such that $V_{gx}
\subseteq U_{hy}$. Then
$V_y \cup V_x \cup V_{gx} \subseteq U_y$, but it is not consistent
with the remark
made above.
\end{pf}

     Thus we have  the following result.

\begin{them}\label{thC}     Let $(X, \rho _0)$ be a  locally
compact and locally  connected metric
space, $G$ a countable group acting discretely on $X$ such that the
quotient space $X/G$ is
metrizable and let us fix a metric $d$ on it.  Then the formula
\label{rho1} $$\rho (x,y) = inf \sum_{i=1}^{k-1} d(px_i, px_{i+1})$$
where the infimum is over all allowable sequences
joining $x$ with $y$
 endows $X$ with a
G-invariant metric $\rho$
topologically equivalent to  $\rho _0$.
\end{them}

The result has several interesting corollaries.

First, it is possible to weaken the assumptions in
Theorem \ref{thC},
since it is enough to assume that X is a metrizable, locally
compact, locally connected space  and
$X/G$ is Hausdorff  to
ensure the existence  of an invariant metric. This follows
from Stone's theorem, since the image of a locally
compact  space under an open map  is
again locally compact.

One has also the following consequence.

\begin{cor}
If a countable group $G$ acts effectively on a compact
connected
manifold $M$ \stt\ $M/G$ is Hausdorff, then $G$ is finite.
\end{cor}

This and the next corollary are related to the classical
results of Newman \cite{N}, Smith \cite{S}, and Montgomery - Zippin
\cite{MZ}, Section 5.5.5.

\begin{cor} Let $M$ be a connected topological
manifold, $G$ a countable discrete
group acting on $M$ \stt\ $M/G$ is Hausdorff. Then

a) (rigidity) any element $g \in G$ whose
set of fixed points has nonempty interior is the identity
transformation,

b) the action
is perfect and all isotropy groups are finite.
\end{cor}

\begin{pf}

a) Let $g$ be any element of $M$ whose set of fixed
points has nonempty interior.
If $g \neq id$ then one can find $x \in \partial(Fix\, g)$
such that
each ball centered at $x$ intersects the interior of
 $Fix\, g$.
It follows from the existence of a relatively compact
slice $U_x$ and the theorems
listed before the corollary that $U_x \subset Fix\, g$,
which is a contradiction.

b) Let $x \in M$ and assume that $G_x$ is infinite.
Then by part a) $G_x$ acts effectively on some
relatively compact slice $U_x$ at $x$.
When we denote $M_g = \{ y \in U_x : gy = y\},$ then by
Baire's theorem there exists $z \in U_x - \bigcup_{g \in G_x} M_g$.
$G_x z $ is infinite in $U_x.$ Thus we obtain a contradiction.

Now we want to show that for any compacts
$A,B \subset X$ the set $\{ g \in G : gA \cap B \neq \emptyset\}$
is finite.
If not, then the existence of slice yields the existence of an
element of $M$ with an infinite stabilizer.
\end{pf}

\begin{cor}
Let $X$ be a compact, locally connected, metrizable $G$-space.
If $G$ is a countable group and the orbit space is Hausdorff,
then there exists an invariant metric on $X.$
\end{cor}

\begin{cor}
Let $X$ be a compact, locally connected
metric $G$-space,
$G$ a countable group. If $X/G$ is Hausdorff,
then any sequence $\{ g_n\}$ in
$G$ contains a subsequence which is convergent in $Homeo\, (X)$
(in the compact-open topology) to a homeomorphism.
\end{cor}

\begin{pf}
We take an invariant metric on $X.$ Using
similar arguments as in the classical Ascoli-Arzela's theorem
we can find a subsequence convergent to
some element of Homeo($X$) in the compact-open topology.
\end{pf}

\begin{cor}
Let $X$ be a locally compact, locally connected, metrizable
and connected
$G$-space. If $G$ is a countable group acting on $X$ such
that $X/G$
is Hausdorff, then either all stabilizers are finite, or
all are infinite.
\end{cor}

\begin{pf}
From the existence of a relatively compact slice at each point
we see that the set of points with infinite stabilizers is open.
Existence of a slice
gives that the set of points with finite stabilizers
is open. So for a
connected space we see that only one type of
stabilizers can appear.
\end{pf}

\newcommand{\gsv}{\mbox{$\tilde{\GG}_s(V)$}}
\newcommand{\gsw}{\mbox{$\tilde{\GG}_s(W)$}}
\newcommand{\gso}{\mbox{$\tilde{\GG}_s(0)$}}
\newcommand{\gsvw}{\mbox{$\tilde{\GG}_s(V\times W)$}}


\begin{thebibliography}{A}


\bibitem [A]{A}
  D.W. Alekseevski,
{\it
On perfect actions of groups (Russian).}
Uspiekhi Math. Nauk {\bf 34} (1979), 219--220


\bibitem [B]{B}
G. Birkhoff,
{\it
A note on topological groups.}
Compositio Math. {\bf 3} (1936),
427-430.


\bibitem [E]{E}
R.  Engelking,
{\bf
General Topology.}  PWN; Warszawa 1977


\bibitem [Ka]{Ka}
S. Kakutani
{\it
\"Uber die Metrization der topologischen Gruppen.}
\newblock
Proc. Imp. Acad. Japan {\bf 12} (1936),  82-84.


\bibitem [Ko]{Ko}
J.L.Koszul,
{\bf
Lectures on Groups of Transformation.}
\newblock
Tate Institute,
Bombay 1965.


\bibitem[MZ]{MZ}
D.Montgomery, L.Zippin
{\bf
Topological Transformation Groups.}
\newblock
Interscience; London, New York, 1995.


\bibitem [N]{N}
M.H. Newman,
{\it
A theorem on periodic transformations of spaces.}
\newblock
Quart. Journal of Math. {\bf 2} (1931), pp. 1-8.


\bibitem[P]{P}
R. S.  Palais,
{\it
On the existence of slices for actions of non-compact
Lie groups.}
\newblock
Ann. Math. {\bf 73} (1961), 295 -- 323.


\bibitem[S]{S}
P.A.Smith,
{\it
Transformations of finite period III.}
\newblock
Ann. Math. {\bf 42} (1941), 446--458.

\end{thebibliography}
\bibliographystyle{amsalpha}

\medskip
\noindent {\bf Mathematical Institute,
Wroc\l aw University,

\noindent pl. Grunwaldzki 2/4,

\noindent 50-384 Wroc\l aw, Poland}

\begin{flushleft}
\tt hajduk@math.uni.wroc.pl

\tt rwalc@math.uni.wroc.pl
\end{flushleft}

\end{document}